\documentclass[10pt]{article}

\usepackage{latexsym}
\usepackage{eucal}
\usepackage{amsxtra}
\usepackage{amsfonts}
\usepackage{selinput}
\usepackage{tikz-cd}
\SelectInputMappings{%
  adieresis={ä},
  eacute={é},
  Lcaron={Ľ},
}

\newtheorem{innercustomthm}{Theorem}
\newenvironment{customthm}[1]
  {\renewcommand\theinnercustomthm{#1}\innercustomthm}
  {\endinnercustomthm}

\newtheorem{theorem}{Theorem}[section]

\newtheorem{corollary}{Corollary}[section]

\newcommand{\Q}{\mbox{$Q^{m}\;$}}

\newcommand{\M}{\mbox{$M\;$}}

\everymath{\displaystyle}

\newenvironment{multi}
{ \begin{equation} \begin{array}{c} }
{ \end{array} \end{equation} }

\makeatletter

\@addtoreset{equation}{section}
\makeatother

\begin{document}

\title{On the structure Lie operator of a real hypersurface in the complex quadric}
\author{Juan de Dios P\'{e}rez, David P\'{e}rez-L\'{o}pez and Young Jin Suh}
\date{}
\maketitle

\begin{abstract}
The almost contact metric structure that we have on a real hypersurface \M in the complex quadric $Q^{m}=SO_{m+2}/SO_mSO_2$ allows us to define, for any nonnull real number $k$, the $k$-th generalized Tanaka-Webster connection on \M, $\hat{\nabla}^{(k)}$. Associated to this connection we have Cho and torsion operators, $F_X^{(k)}$ and $T_X^{(k)}$, respectively, for any vector field $X$ tangent to $M$. From them and for any symmetric operator $B$ on $M$ we can consider two tensor fields of type (1,2) on $M$ that we will denote by $B_F^{(k)}$ and $B_T^{(k)}$, respectively. We will classify real hypersurfaces $M$ in $Q^m$ for which any of those tensors identically vanishes, in the particular case of $B$ being the structure Lie operator $L_{\xi}$ on $M$.
\end{abstract}



\section{Introduction.}

Let $M$ be a real hypersurface of the complex quadric $Q^m=SO_{m+2}/SO_mSO_2$ endowed with the Kählerian structure $(J,G)$. Let $N$ be a unit local normal vector field on $M$. $(J,G)$ induces on $M$ an almost contact metric structure $(\phi ,\eta ,g,\xi)$, see \cite{B}. Then, if $\nabla$ denotes its Levi-Civita connection, for any nonnull real number $k$, we can define on $M$ the so-called $k$-th generalized Tanaka-Webster connection $\hat{\nabla}^{(k)}$ by

\[  \hat{\nabla}_X^{(k)}Y=\nabla_XY+g(\phi SX,Y)\xi -\eta(Y)\phi SX-k\eta(X)\phi Y   \]

\noindent for any $X,Y$ tangent to \M (see \cite{CHO1999}, \cite{CHO2008}, \cite{T}, \cite{TAN} and \cite{W}), where $S$ denotes the shape operator on $M$ associated to $N$.

We will write $F_X^{(k)}Y=g(\phi SX,Y)\xi -\eta(Y)\phi SX-k\eta(X)\phi Y$, for any $X,Y$ tangent to \M and call it the $k$-th Cho operator on $M$ associated to $X$. The torsion of the $k$-th generalized Tanaka-Webster connection is given by $T^{(k)}(X,Y)=F_X^{(k)}Y-F_Y^{(k)}X$, for any $X,Y$ tangent to $M$. We will call $k$-th torsion operator associated to the vector field $X$ to the operator given by $T_X^{(k)}Y=T^{(k)}(X,Y)$ for any $Y$ tangent to $M$.                                                                                                                         

If $\cal L$ denotes the Lie derivative on $M$, we have that ${\cal L}_XY=\nabla_XY-\nabla_YX$, for any $X,Y$ tangent to $M$. We can also consider on $M$ a differential operator of first order associated to $\hat{\nabla}^{(k)}$, ${\cal L}^{(k)}$, given by ${\cal L}_X^{(k)}Y=\hat{\nabla}_X^{(k)}Y-\hat{\nabla}_Y^{(k)}X$, for any $X,Y$ tangent to $M$. It is easy to see that ${\cal L}_X^{(k)}Y={\cal L}_XY+T_X^{(k)}Y$, for any $X,Y$ tangent to $M$.

Let $B$ be a symmetric operator on $M$. Then we can consider $((\hat{\nabla}_X^{(k)}-\nabla_X)B)Y$ for any $X$ and $Y$ tangent to $M$. It is easy to see that $((\hat{\nabla}_X^{(k)}-\nabla_X)B)Y=F_X^{(k)}BY-BF_X^{(k)}Y$. Therefore, from $B$, we can define on $M$ the tensor field of type (1,2) $B_F^{(k)}$ that is given by

\[    B_F^{(k)}(X,Y)=\lbrack F_X^{(k)},B \rbrack Y=F_X^{(k)}BY-BF_X^{(k)}Y    \]

\noindent for any $X,Y$ tangent to $M$. Thus the vanishing of such a tensor field is equivalent to the commutativity of the operators $F_X^{(k)}$ and $B$ for any $X$ tangent to $M$. We have studied this problem in the cases of $B=S$, the shape operator of $M$, in \cite{PS2018} and $B=\bar{R}_N$, the normal Jacobi operator of $M$, in \cite{PS2019}.

We can also consider $(({\cal L}_X^{(k)}-{\cal L}_X)B)Y$ for any $X,Y$ tangent to $M$. As above, it is easy to see that $(({\cal L}_X^{(k)}-{\cal L}_X)B)Y=T_X^{(k)}BY-BT_X^{(k)}Y$, and we can define another tensor field of type (1,2) on $M$ associated to $B$, $B_T^{(k)}$, by

\[     B_T^{(k)}(X,Y)=\lbrack T_X^{(k)},B \rbrack Y=T_X^{(k)}BY-BT_X^{(k)}Y   \]

\noindent for any $X,Y$ tangent to $M$. Such a tensor field vanishes if and only if, for any $X$ tangent to $M$, the operators $T_X^{(k)}$ and $B$ commute. This problem has been solved in \cite{PS2018} for $B=S$ and in \cite{P} for $B=\bar{R}_N$.

Let us define the structure Lie operator on $M$, $L_{\xi}$, by $g(L_{\xi}X,Y)=({\cal L}_{\xi}g)(X,Y)=g(\nabla_X\xi ,Y)+g(X,\nabla_Y,\xi)=g(\phi SX,Y)+g(X,\phi SY)=g((\phi S-S\phi)X,Y)$, for any $X,Y$ tangent to $M$. Therefore we have

\[   L_{\xi}X=(\phi S-S\phi)X              \]

\noindent for any $X$ tangent to $M$. Berndt and Suh, \cite{BS2013}, classified real hypersurfaces in $Q^m$ with isometric Reeb flow, proving that such a condition is equivalent to the vanishing of the structure Lie operator. They obtained the following

\begin{customthm}{A}
Let $M$ be a real hypersurface in the complex quadric $Q^m$, $m \geq 3$. The Reeb flow on $M$ is isometric if and only if $m$ is even, say $m=2n$, and $M$ is an open part of a tube around a totally geodesic $\mathbb{C} P^n \subset Q^m$.
\end{customthm}

A real hypersurface $M$ in $Q^m$ is Hopf if the Reeb vector field $\xi$ is principal, that is, $S\xi =\alpha\xi$, and $\alpha$ is called the Reeb function of $M$. The tubes appearing in Theorem 1.1 are examples of Hopf real hypersurfaces in $Q^m$. We will denote by $\mathcal{C}$ the maximal holomorphic distribution on $M$, given by ${\mathcal C}=\{ X \in TM |$ $g(X,\xi)=0\}$.

Recently, in \cite{KLPS}, a family of non Hopf real hypersurfaces of $Q^m$, called ruled real hypersurfaces, was introduced and some examples were provided. $M$ is ruled in $Q^m$ if the distribution $\mathcal C$ is integrable and $M$ is foliated by totally geodesic complex hyperplanes $Q^{m-1}$ in $Q^m$. In such a case, the shape operator $S$ of $M$ in $Q^n$ is given by $S\xi =\alpha\xi +\beta U$, $SU=\beta\xi$ and $SX=0$ for any vector field $X \perp Span\{\xi ,U\}$, where $U$ is a unit vector field in $\mathcal C$, $\alpha$ and $\beta$ are functions on $M$ and $\beta$ does not vanish.

This paper is devoted to study the vanishing of the tensors $L_{\xi_F}^{(k)}$ and $L_{\xi_T}^{(k)}$. We will prove the following

\begin{theorem}
Let \M be a real hypersurface in \Q, $m \geq 3$, and $k$ a nonnull real number. Then $L_{\xi_F}^{(k)}(\xi ,X)=0$, for any $X$ tangent to $M$, if and only if $M$ is an open part of a tube of radius $r$, $0 < r <\frac{\pi}{2}$, around a totally geodesic $\mathbb{C}P^n \subset Q^m$, $m=2n$. 
\end{theorem}

As $TM=\mathbb{R}\xi \oplus \mathcal{C}$, we can prove also the 

\begin{theorem}
Let \M be a real hypersurface of \Q, $m \geq 3$, and $k$ a nonnull real number. Then $L_{\xi_F}^{(k)}(X ,Y)=0$ for any $X \in \mathcal{C}$, $Y$ tangent to $M$ if and only if \M is an open part of either a tube of radius $r$, $0 < r <\frac{\pi}{2}$, around a totally geodesic $\mathbb{C}P^n \subset Q^m$, $m=2n$, or a ruled real hypersurface.
\end{theorem}

From both Theorems we obtain the

\begin{corollary}
Let $M$ be a real hypersurface in \Q, $m \geq 3$, and $k$ a nonnull real number. Then $L_{\xi_F}^{(k)}\equiv 0$ if and only if $M$ is an open part of a tube of radius $r$, $0 < r <\frac{\pi}{2}$, around a totally geodesic $\mathbb{C}P^n \subset Q^m$, $m=2n$.
\end{corollary}

In the case of $L_{\xi_T}^{(k)}$ we have, bearing in mind again that $TM=\mathbb{R}\xi \oplus \mathcal{C}$, we have

\begin{theorem}
Let \M be a real hypersurface in \Q, $m \geq 3$, and $k$ and nonnull real number. Then $L_{\xi_T}^{(k)}(X,Y)=$0 for any $X \in \mathcal{C}$, $Y$ tangent to $M$, if and only if $M$ is an open part of a tube of radius $r$, $0< r <\frac{\pi}{2}$, around $\mathbb{C}P^n$, $m=2n$.
\end{theorem}

On the other hand

\begin{theorem}
There does not exist any non Hopf real hypersurface $M$ in \Q, $m \geq 3$, such that  $L_{\xi_T}^{(k)}(\xi ,X)=0$ for any $X$ tangent to \M and any nonnull real number $k$.
\end{theorem}

As a consequence of these results we conclude with the

\begin{corollary}
Let  \M be a real hypersurface in \Q, $m \geq 3$, and $k$ a nonnull real number. Then $L_{\xi_T}^{(k)} \equiv 0$ if and only if $M$ is an open part of a tube of radius $r$, $0< r <\frac{\pi}{2}$, around $\mathbb{C}P^n$, $m=2n$. 
\end{corollary}

\section{Preliminaries.}

For more details about this section we refer to  \cite{BS2012}, \cite{BS2013}, \cite{BSB}, \cite{K}, \cite{KO} and \cite{R}. The complex quadric $Q^m$ is the complex hypersurface in the complex projective space $\mathbb{C}P^{m+1}$ defined by the equation $z_1^2 + \cdots + z_{m+2}^2 = 0$, where $z_1,\ldots,z_{m+2}$ are homogeneous coordinates on $\mathbb{C}P^{m+1}$. The Kähler structure on $\mathbb{C}P^{m+1}$ induces canonically a Kähler structure $(J,G)$ on the complex quadric, where $G$ is a Riemannian metric induced on $Q^m$ by the Fubini-Study metric of $\mathbb{C}P^{m+1}$.

A point $[z]$ in $\mathbb{C}P^{m+1}$ is given by $[z]=\{ \lambda z |$ $z \in \mathbb{C} \}$, where $z$ is a nonzero vector in $\mathbb{C}^{m+2}$. If $\bar{z}$ is the conjugated of $z$ , $T_{[z]}Q^m$, for any $[z] \in Q^m$, can be identified canonically with $\mathbb{C}^{m+2} \ominus ([z] \oplus [\bar{z}])$, the orthogonal complement of $[z] \oplus [\bar{z}]$ in $\mathbb{C}^{m+2}$. $\bar{z}$ is a unit normal vector of $Q^m$ in $\mathbb{C}P^{m+1}$ at the point $[z]$. The shape operator $A_{\bar{z}}$ of $Q^m$ with respect to the unit normal vector $\bar{z}$ is given by

\[   A_{\bar{z}}w=\bar{w}         \]

\noindent for any $w \in T_{[z]}Q^m$. Then $A_{\bar{z}}$ defines a real structure on the complex vector space $T_{[z]}Q^m$. The set of all shape operators $A_{\lambda \bar{z}}$ defines a parallel subbundle of rank 2, $\mathfrak A$, of $End(TQ^m)$, consisting of all the real structures on the tangent space of $Q^m$. For any $A \in \mathfrak{A}$, $A^2=I$, the identity endomorphism, and $AJ=-JA$.

A nonzero tangent vector $W \in T_{[z]}Q^m$ is called singular if it is tangent to more than one maximal flat in $Q^m$. There are two types of singular tangent vectors in $Q^m$:
\begin{enumerate}
\item $W$ is $\mathfrak A$-principal if there exists a real structure $A \in \mathfrak{A}$ such that $AW=W$.

\item $W$ is $\mathfrak{A}$-isotropic if there exist $A \in \mathfrak{A}$ and two orthonormal vectors $X,Y$ such that $AX=X$, $AY=Y$ and $W/\| W\| =(X+JY)/\sqrt{2}$.
\end{enumerate}

If $M$ is a real hypersurface in \Q with unit local normal vector field $N$ and Levi-Civita connection $\nabla$, let $S$ denote the shape operator on $M$ with respect to $N$. For any vector field $X$ tangent to \M we write $JX=\phi X+\eta(X)N$. The tangential component of $JX$, $\phi X$ defines on $M$ a skew-symmetric tensor field of type (1,1)  called the structure tensor. The vector field  $\xi=-JN$ is called the Reeb vector field (or structure vector field) of \M. Thus we naturally obtain that the 1-form $\eta$ is given by $\eta(X)=g(X,\xi)$ for any vector field $X$ tangent to \M. Therefore on $M$ we have an almost contact metric structure ($\phi,\xi,\eta,g$), where $g$ is the restriction to $M$ of the metric $G$. Thus, see \cite{B}, we have the following relations

\begin{equation} \label{2.1}
\phi^2X=-X+\eta(X)\xi,  \quad  \eta(\xi)=1,  \quad   g(\phi X,\phi Y)=g(X,Y)-\eta(X)\eta(Y)
\end{equation}
\noindent for any $X,Y$ tangent to $M$. From (\ref{2.1}) we also have
\begin{equation}\label{2.2}
\phi\xi =0,   \quad   \nabla_X\xi =\phi SX
\end{equation}

\noindent for any $X$ tangent to $M$, where we have applied that $J$ is parallel respect to the Levi-Civita connection on $Q^m$

The tangent bundle $TM$ of $M$ splits orthogonally into
\[TM=\mathcal{C}\oplus\mathcal{F},\]
where $\mathcal{C}=ker(\eta)=\{ X \in TM | g(X,\xi)=0 \}$ is the maximal complex (holomorphic) subbundle of $TM$ and $\mathcal{F}=\mathbb{R}\xi$. If we restrict $\phi$ to $\mathcal{C}$ it coincides with the complex structure $J$.

For a tube M of radius r, $0<r<\pi/2$, around the totally geodesic $\mathbb{C}P^{n}$ in \Q, $m=2n$, its normal bundle consists of $\mathfrak{A}$-isotropic singular tangent vectors  of \Q, see \cite{BS2013}.

In the case of a ruled real hypersurface of $Q^m$, $N$ is $\mathfrak{A}$-principal, see \cite{KLPS}.

\section{Proofs of Theorems 1.1 and 1.2.}

Suppose that $M$ satisfies $L_{\xi_F}^{(k)}(\xi ,X)=0$ for any $X$ tangent to $M$. This means that $F_{\xi}^{(k)}(\phi S-S\phi)X-(\phi S-S\phi)F_{\xi}^{(k)}X=0$ for any $X$ tangent to $M$.  Therefore,

\begin{multi}\label{3.1}
g(\phi S\xi ,(\phi S-S\phi)X)\xi -\eta((\phi S-S\phi)X)\phi S\xi -k\phi(\phi S-S\phi)X-g(\phi S\xi ,X)\phi S\xi   \\
+\eta(X)(\phi S-S\phi)\phi S\xi +k(\phi S-S\phi)\phi X=0
\end{multi}

\noindent for any $X$ tangent to $M$. Suppose first that \M is Hopf with Reeb function $\alpha$. In this case (\ref{3.1}) yields

\begin{multi}\label{3.2}
-k\phi(\phi S-S\phi)X+k(\phi S-S\phi)\phi X=0
\end{multi}

\noindent for any $X$ tangent to $M$. As $k \neq 0$, from (\ref{3.2}) we get $2SX+2\phi S\phi X=0$ for any $X \in \mathcal{C}$. If we apply $\phi$  to this equation we obtain $S\phi X=\phi SX$ for any $X \in \mathcal{C}$. This yields $S\phi =\phi S$ and the results follows from Theorem A.

If $M$ is non Hopf, we write $S\xi=\alpha\xi +\beta U$, where $U$ is a unit vector field in $\mathcal C$, $\alpha$ and $\beta$ are functions on $M$ and $\beta$ does not vanish at least on a neighbourhood of a point $p \in M$. The following calculations will be made on such a neighbourhood. From (\ref{3.1}) we have

\begin{multi}\label{3.3}
\beta g(SU,X)\xi -\beta g(S\phi U,\phi X)\xi-2\beta^2g(\phi U,X)\phi U-k\phi(\phi S-S\phi)X  \\
+\beta\eta(X)(\phi S-S\phi)\phi U+k(\phi S-S\phi)\phi X=0
\end{multi}

\noindent for any $X$ tangent to $M$. Taking $X=\xi$ in (\ref{3.3}) it follows that $\beta^2\xi +k\beta U+\beta(\phi S-S\phi)\phi U=0$. Its scalar product with $\xi$ yields $2\beta^2=0$, which is impossible and finishes the proof of Theorem 1.1 

If we now suppose that $L_{\xi_F}^{(k)}(X,Y)=0$, for any $X \in \mathcal{C}$, $Y$ tangent to $M$, we obtain

\begin{multi}\label{3.4}
g(\phi SX,\phi SY)\xi -g(\phi SX,S\phi Y)\xi -g((\phi S-S\phi)Y,\xi)\phi SX  \\
-g(\phi SX,Y)\phi S\xi +\eta(Y)(\phi S-S\phi)\phi SX=0  \\
\end{multi}

\noindent for any $X \in \mathcal{C}$, $Y$ tangent to $M$. First suppose that $M$ is Hopf with Reeb function $\alpha$. Then (\ref{3.4}) becomes

\begin{multi}\label{3.5}
g(SX,SY)\xi -g(\phi SX,S\phi Y)\xi +\eta(Y)(\phi S-S\phi)\phi SX=0
\end{multi}

\noindent for any $X \in \mathcal{C}$, $Y$ tangent to $M$. Take $Y \in \mathcal{C}$ in (\ref{3.5}). Then $g(S^2X,Y)\xi +g(\phi S\phi SX,Y)\xi =0$ for any $X,Y \in \mathcal{C}$. This means that  $S^2X+\phi S\phi SX=0$ for any $X \in \mathcal{C}$. Apply $\phi$ to obtain $\phi S^2X=S\phi SX$, for any $X \in \mathcal{C}$. Suppose that $X \in \mathcal{C}$ satisfies $SX=\lambda X$. Then we obtain $\lambda^2\phi X=\lambda S\phi X$. Therefore, either $\lambda =0$ or $\lambda\neq 0$ and $S\phi X=\lambda \phi X$. This means that for any nonnull eigenvalue of $S$ in $\mathcal{C}$ the corresponding eigenspace is $\phi$-invariant. Thus also the eigenspace corresponding to the possible eigenvalue $0$ is also $\phi$-invariant and $S\phi =\phi S$. We conclude as in Theorem 1.1.

If we suppose now that $M$ is non Hopf, we write, as above, $S\xi =\alpha\xi +\beta U$, with the same conditions. From (\ref{3.4}) we have

\begin{multi}\label{3.6}
g(\phi SX,\phi SY)\xi -g(\phi SX,S\phi Y)\xi -\beta g(Y,\phi U)S\phi X     \\
-\beta g(\phi SX,Y)\phi U+\eta(Y)(\phi S-S\phi)\phi SX=0
\end{multi}

\noindent for any $X \in \mathcal{C}$, $Y$ tangent to $M$. If we take $Y=\xi$ in (\ref{3.6}) we get $\beta g(SU,X)\xi +(\phi S-S\phi)\phi SX=0$, for any $X \in \mathcal{C}$. Its scalar product with $\xi$ yields $2\beta g(SU,X)=0$, and as $\beta \neq 0$, we obtain $g(SU,X)=0$, for any $X \in \mathcal{C}$. That is,

\begin{multi}\label{3.7}
SU=\beta\xi
\end{multi}

If we take $Y=\phi U$ in (\ref{3.6}) it follows $g(\phi SX,\phi S\phi U)\xi -\beta\phi SX=0$, for any $X \in \mathcal{C}$. Its scalar product with $U$ implies $\beta g(S\phi U,X)=0$ for any $X \in \mathcal{C}$. $\beta$ being non vanishing implies

\begin{multi}\label{3.8}
S\phi U=0.
\end{multi}

We will denote $\mathcal{C}_U=\{ X \in \mathcal{C} |$ $g(X,U)=g(X,\phi U)=0 \}$. From (\ref{3.7}) and (\ref{3.8}) we know that $\mathcal{C}_U$ is $S$-invariant. Take $X \in \mathcal{C}_U$, $Y=\phi U$ in (\ref{3.6}). This gives $-\beta\phi SX=0$. Therefore, $\phi SX=0$ for any $X \in \mathcal{C}_U$, and applying $\phi$ we get

\begin{multi}\label{3.9}
SX=0
\end{multi}

\noindent for any $X \in {\mathcal{C}}_U$. From (\ref{3.7}), (\ref{3.8}) and (\ref{3.9}) we obtain that $M$ is ruled. The converse is trivial and we have finished the proof of Theorem 1.2.

\section{Proofs of Theorems 1.3 and 1.4.}

If we suppose that $T_X^{(k)}L_{\xi}Y-L_{\xi}T_X^{(k)}Y=0$, for any $X  \in \mathcal{C}$, $Y$ tangent to $M$, we should have $F_X^{(k)}(\phi S-S\phi)Y-F_{(\phi S-S\phi) Y}^{(k)}X-(\phi S-S\phi)F_X^{(k)}Y+(\phi S-S\phi)F_Y^{(k)}X=0$, for any $X \in \mathcal{C}$, $Y$ tangent to $M$. Therefore,

\begin{multi}\label{4.1}
g(\phi SX,(\phi S-S\phi)Y)\xi -\eta((\phi S-S\phi)Y)\phi SX-g(\phi S(\phi S-S\phi)Y,X)\xi  \\
+k\eta((\phi S-S\phi)Y)\phi X-g(\phi SX,Y)\phi S\xi +\eta(Y)(\phi S-S\phi)\phi SX    \\
+g(\phi SY,X)\phi S\xi -\eta(Y)(\phi S-S\phi)\phi X =0  \\
\end{multi}

\noindent for any $X \in \mathcal{C}$, $Y$ tangent to $M$. Suppose first that $M$ is Hopf and write $S\xi =\alpha\xi$. Then (\ref{4.1}) gives

\begin{multi}\label{4.2}
g(SX,SY)\xi -g(\phi SX,S\phi Y)\xi -g(\phi S(\phi S-S\phi)Y,X)\xi   \\
+\eta(Y)(\phi S-S\phi)\phi SX-k\eta(Y)(\phi S-S\phi)\phi X=0
\end{multi}

\noindent for any $X \in \mathcal{C}$, $Y$ tangent to $M$. If we take $Y=\xi$ in (\ref{4.2}) we obtain $(\phi S-S\phi)\phi SX-k(\phi S-S\phi)\phi X=0$, for any $X \in \mathcal{C}$. That is, $\phi S\phi SX+S^2X-k\phi S\phi X-kSX=0$, for any $X \in \mathcal{C}$. If we suppose that $X \in \mathcal{C}$ satisfies $SX=\lambda X$, we have $\lambda\phi S\phi X+\lambda^2X-k\phi S\phi X-k\lambda X=0$. Applying $\phi$ we get $-\lambda S\phi X+\lambda^2\phi X+kS\phi X-k\lambda\phi X=0$. Then, $(k-\lambda)S\phi X=\lambda(k-\lambda)\phi X$. Thus either $\lambda =k$ or, if $\lambda \neq k$, $S\phi X=\lambda \phi X$. We have found that for any eigenvalue $\lambda \neq k$ in $\mathcal{C}$, the corresponding eigenspace is $\phi$-invariant. This yields that also the eigenspace corresponding to $k$ is $\phi$-invariant. Therefore $\phi S=S\phi$ and $M$ is an open part of a tube of radius $r$, $0 < r < \frac{\pi}{2}$ around a totally geodesic $\mathbb{C}P^n \subset Q^m$, with $m=2n$.

Let us suppose now that $M$ is non Hopf and write again $S\xi =\alpha\xi +\beta U$ as in previous section. From (\ref{4.1}) we get

\begin{multi}\label{4.3}
g(\phi SX,\phi SY)\xi -g(\phi SX,S\phi Y)\xi -\beta g(Y,\phi U)\phi SX+g((\phi S-S\phi)Y,S\phi X)\xi    \\
+k\beta g(Y,\phi U)\phi X-\beta g(\phi SX,Y)\phi U +\eta(Y)(\phi S-S\phi)\phi SX +\beta g(\phi SY,X)\phi U   \\
-k\eta(Y)(\phi S-S\phi)\phi X=0    \\
\end{multi}

\noindent for any $X \in \mathcal{C}$, $Y$ tangent to $M$. If we take $X=\xi$ in (\ref{4.3}) it follows

\begin{multi}\label{4.4}
\beta g(SU,X)\xi +\beta g(S\phi U,\phi X)\xi +(\phi S-S\phi)\phi SX+\beta^2g(\phi U,X)\phi U  \\
-k(\phi S-S\phi)\phi X=0
\end{multi}

\noindent for any $X \in \mathcal{C}$. The scalar product of (\ref{4.4}) and $\xi$ gives $\beta g(SU,X)+\beta g(S\phi U,\phi X)+g(\phi SX,\phi S\xi)-kg(\phi X,\phi S\xi)=0$ for any $X \in \mathcal{C}$. Therefore, $2g(SU,X)+g(S\phi U,\phi X)-kg(U,X)=0$ for any $X \in \mathcal{C}$. For $X=\phi U$ we obtain

\begin{multi}\label{4.5}
g(SU,\phi U)=0
\end{multi}

\noindent and for $X=U$ it yields

\begin{multi}\label{4.6}
2g(SU,U)+g(S\phi U,\phi U)=k.
\end{multi}

Taking $Y \in \mathcal{C}$ in (\ref{4.3}) we have

\begin{multi}\label{4.7}
g(\phi SX,\phi SY)\xi -g(\phi SX,S\phi Y)\xi -\beta g(Y,\phi U)\phi SX+g((\phi S-S\phi)Y,S\phi X)\xi   \\
+k\beta g(Y,\phi U)\phi X-\beta g(\phi SX,Y)\phi U+\beta g(\phi SY,X)\phi U=0
\end{multi}

\noindent for any $X,Y \in \mathcal{C}$. Its scalar product with $\phi U$ gives

\begin{multi}\label{4.8}
-g(Y,\phi U)g(SX,U)+kg(Y,\phi U)g(X,U)-g(\phi SX,Y)+g(\phi SY,X)=0
\end{multi}

\noindent for any $X,Y \in \mathcal{C}$. If in (\ref{4.8}) we take $Y=U$ we obtain $-g(\phi SX,U)+g(S\phi U,X)=0$ for any $X \in \mathcal{C}$. Thus

\begin{multi}\label{4.9}
S\phi U+\phi SU=0.
\end{multi}

Taking $Y=\phi U$ in (\ref{4.8}) it follows $-2g(SU,X)+kg(U,X)+g(\phi S\phi U,X)=0$ for any $X \in \mathcal{C}$. Then $-2SU+kU+\phi S\phi U$ has no component in $\mathcal{C}$. That means that $-2SU+kU+\phi S\phi U=g(-2SU+kU+\phi S\phi U,\xi)\xi =-2\beta\xi$. On the other hand, from (\ref{4.9}), $-2\beta\xi =-2SU+kU-\phi^2SU=-SU+kU-\beta\xi$. That is,

\begin{multi}\label{4.10}
SU=\beta\xi +kU.
\end{multi}

From (\ref{4.10}) and (\ref{4.9}),

\begin{multi}\label{4.11}
S\phi U=-k\phi U.
\end{multi}

Therefore, $\mathcal{C}_U$ is $S$-invariant. Taking $X \in \mathcal{C}_U$ in (\ref{4.7}) we get

\begin{multi}\label{4.12}
g(S^2Y,X)\xi +g(S\phi S\phi Y,X)\xi -g(\phi S(\phi S-S\phi)Y,X)\xi    \\
+\beta g(S\phi Y,X)\phi U+\beta g(\phi SY,X)\phi U=0
\end{multi}

\noindent for any $X \in \mathcal{C}$, $Y \in \mathcal{C}_U$. Its scalar product with $\phi U$ yields $g(S\phi Y,X)+g(\phi SY,X)=0$, for any $Y \in \mathcal{C}_U$, $X \in \mathcal{C}$. Then

\begin{multi}\label{4.13}
S\phi Y+\phi SY=0
\end{multi}

\noindent for any $Y \in \mathcal{C}_U$.

Taking $Y=\phi U$ in (\ref{4.3}) we obtain $g(SX,S\phi U)\xi +g(SU,\phi SX)\xi -\beta\phi SX+g((\phi S-S\phi)\phi U,S\phi X)\xi +k\beta\phi X-\beta g(SX,U)\phi U+\beta g(\phi S\phi U,X)\phi U=0$, for any $X \in \mathcal{C}$. If $X \in \mathcal{C}_U$ this yields $-\beta\phi SX+k\beta\phi X=0$. As we suppose $\beta \neq 0$, $SX=kX$, for any $X \in \mathcal{C}_U$, but then (\ref{4.13}) yields $2k=0$, which is impossible, finishing the proof of Theorem 1.3.

If now we suppose that $M$ is non Hopf and $L_{\xi_T}^{(k)}(\xi ,X)=0$ for any $X$ tangent to $M$, we obtain

\begin{multi}\label{4.14}
\beta g(\phi U,(\phi S-S\phi)X)\xi -2\beta^2(X,\phi U)\phi U-k\phi (\phi S-S\phi)X     \\
+\phi S(\phi S-S\phi)X+\beta\eta(X)(\phi S-S\phi)\phi U+k(\phi S-S\phi)\phi X      \\
-(\phi S-S\phi)\phi SX=0
\end{multi}

\noindent for any $X$ tangent to $M$. Taking $X=\xi$ in (\ref{4.14}) we get $\beta^2\xi-k\phi^2S\xi +\phi S\phi S\xi +\beta(\phi S-S\phi)\phi U-\beta(\phi S-S\phi)\phi U=0$. Therefore, $\beta^2\xi +k\beta U+\beta\phi S\phi U=0$. Its scalar product with $\xi$ yields $\beta =0$, a contradiction that finishes the proof of Theorem 1.4.

In order to prove Corollary 1.2, if $M$ satisfies $L_{\xi_T}^{(k)} \equiv 0$, from Theorems 1.3 and 1.4, $M$ must be Hopf. Now from Theorem 1.3 $\phi S=S\phi$, that is, $L_{\xi} \equiv 0$. Therefore $L_{\xi_T}^{(k)}(\xi ,X)=0$, for any nonnull real number $k$ and any $X$ tangent to $M$ and the proof finishes.

\section*{Acknowledgements}
First author is supported by Projects  PID 2020-11 6126GB-I00 from MICINN and PY20-01391 from Junta de Andaluc\'{i}a. Third author by grant Proj. No. NRF-2018-R1D1A1B-05040381 from National Research Foundation of Korea.

\begin{trivlist}

\item Juan de Dios P\'erez: jdperez@ugr.es  \\
Departamento de Geometr\'{i}a y Topolog\'{i}a  and IMAG (Instituto de Matemáticas)\\
Universidad de Granada \\
18071 Granada  \\
Spain  \\

\item David P\'erez-L\'opez: davidpl109@correo.ugr.es \\
18003 Granada  \\
Spain \\

\item Young Jin Suh: yjsuh@knu.ac.kr  \\
Kyungpook National University  \\
College of Natural Sciences  \\
Department of Mathematics  \\
and Research Institute of Real and Complex Manifolds \\
Daegu 41566, Republic of Korea
\end{trivlist}

\end{document}